\documentclass[12pt]{article}
\usepackage{latexsym,amsfonts,amsmath,theorem,amssymb}
\usepackage{graphicx}
\usepackage{color}
\usepackage{diagbox}

\usepackage{enumerate}
\usepackage{hyperref}

\usepackage{algpseudocode,algorithm} % new definition
\algdef{SE}[DOWHILE]{Do}{doWhile}{\algorithmicdo}[1]{\algorithmicwhile\ #1}

\newtheorem{theorem}{Theorem}

\newtheorem{definition}[theorem]{Definition}
{\theoremstyle{plain} \theorembodyfont{\rmfamily}
}
{\theoremstyle{plain} \theorembodyfont{\rmfamily}
\newtheorem{example}[theorem]{Example}}

\allowdisplaybreaks

\newcommand{\bsgamma}{{\boldsymbol{\gamma}}}

\newcommand{\bst}{{\boldsymbol{t}}}

\newcommand{\bsx}{{\boldsymbol{x}}}

\newcommand{\bszero}{{\boldsymbol{0}}}

\newcommand{\rd}{\,\mathrm{d}}

\newcommand{\bbR}{\mathbb{R}}

\newcommand{\bbN}{\mathbb{N}}

\newcommand{\calD}{\mathcal{D}}

\newcommand{\calF}{\mathcal{F}}

\newcommand{\calS}{\mathcal{S}}

\newcommand{\mask}[1]{}

\newcommand{\e}{{\varepsilon}}
\newcommand{\setu}{{\mathfrak{u}}}
\newcommand{\setv}{{\mathfrak{v}}}

\newcommand{\setU}{{\mathfrak{U}}}

\title{Small Superposition Dimension\\
  and Active Set Construction\\
  for Multivariate Integration
  \\ Under Modest Error Demand
  }
\author{A.~D.~Gilbert and G.~W.~Wasilkowski}
\date{\today}

\begin{document}
\maketitle
{\bf Abstract:\ } Constructing {\em active sets} is a key part of the 
Multivariate Decomposition Method. An algorithm for constructing 
optimal or quasi-optimal active sets is proposed in the paper. By numerical 
experiments, it is shown that the new method can provide sets that 
are significantly smaller than the sets constructed by the
already existing method. The experiments also show that the 
{\em superposition dimension} could surprisingly be very small, at most 3, 
when the error demand is not smaller than $10^{-3}$ and the weights decay sufficiently
fast.

\section{Introduction}
In this short paper, we consider approximating integrals with
infinitely many variables. 
We focus on approximations with a modest error demand, aiming at 
problems in, e.g., {\em Mathematical Finance} and
{\em Uncertainty Quantification}, where the underlying stochastic 
process is not known and hence only rough approximations are needed.
In our tests we use $\e=10^{-n}$ for $n=1,2,3$ as the error demands. 

The functions to be integrated belong to
$\{\gamma_\setu\}_{\setu\subset\bbN_+}$-weighted tensor product 
Banach spaces $\calF_\bsgamma$
which allow for the decomposition 
\[
  f(\bsx)\,=\,\sum_{\setu\subset\bbN_+,|\setu|<\infty}f_\setu(\bsx).
\]
Here the summation is with respect to finite subsets $\setu$
of positive
integers and each $f_\setu$ depends only on the variables $x_j$ with
$j\in\setu$. We also assume that the weights $\gamma_\setu$ have a
product form. 

Integrals of such functions can be approximated by the
{\em Multivariate Decomposition Method}, which is a refined version of the
{\em Changing Dimension Algorithm} introduced in \cite{KSWW10}. An essential
part of those methods is the construction of an {\em active set}
$\setU(\e)$ of subsets $\setu$ such that the integral of
$\sum_{\setu\notin\setU(\e)}f_\setu$ can be neglected since it is
bounded from above by
\[
\e\,\left\|\sum_{\setu\notin\setU(\e)}f_\setu
\right\|_{\calF_\bsgamma} \quad\mbox{for all\ }f\,\in\,\calF_\gamma.
\]
In other words, it is enough to approximate integrals of the partial
sum
\[
\sum_{\setu\in\setU(\e)}f_\setu.
\]
We would like to construct possibly small
active sets and such that the largest cardinality among its 
elements $\setu$,
\[
   d(\setU(\e))\,:=\,\max_{\setu\in\setU(\e)}|\setu|,
\]
is also small. This is because the partial sum
$\sum_{\setu\in\setU(\e)}f_\setu$, that can be considered
instead of the infinite sum
$\sum_{\setu\subset\bbN_+,|\setu|<\infty}$, 
has a small 
number $|\setU(\e)|$ of functions $f_\setu$, each depending on no more
than $d(\setU(\e))$ variables. 

A specific construction of such sets (denoted by
$\setU^{\rm PW}(\e)$) was proposed in \cite{PlaWas10} and it was
shown there that the largest cardinality among all 
$\setu\in\setU^{\rm PW}(\e)$ grows very slowly with decreasing $\e$,
\[
d\left(\setU^{\rm PW}(\e)\right)\,=\, 
O\left(\frac{\ln(1/\e)}{\ln(\ln(1/\e))}\right)
\quad\mbox{as\ }\e\,\to\,0.
\]
Moreover the 
size $|\setU^{\rm PW}(\e)|$ grows polynomially in $1/\e$. 
However, the asymptotic constants in the big-$O$ notation 
were not investigated and, as we shall see, they could
be very large. 

This is why in this paper we consider constructing possibly smallest 
active sets denoted by $\setU^{\rm opt}(\e)$.
As we will show by examples, the difference
between the size of $\setU^{\rm opt}(\e)$ and $\setU^{\rm PW}(\e)$
could be very large. 
We also provide a construction of {\em quasi-optimal}
sets, denoted by $\setU^{\rm q-opt}(\e)$, which sometimes are
only slightly larger than the
optimal $\setU^{\rm opt}(\e)$; however, their construction is less expensive.

We are also interested in active sets with the smallest
$d(\setU(\e))$. This leads to the following concept of 
{\em $\e$-superposition dimension} (or {\em superposition dimension}
for short) defined by 
\[
d^{\rm sup}(\e)\,:=\,\min\{d(\setU(\e))\ :\ \setU(\e)\mbox{ is an active set}\}.
\]
Since the optimal active sets in our experiments have very small
$d(\setU^{\rm opt}(\e))$, this implies that the superposition dimension
is also small.

Note that our concept of the superposition dimension depends on the
integration problem as well as the error demand $\e$.  
Hence it is in the same
spirit as the definition of {\em truncation dimension}
introduced recently in
\cite{KPW15}. They are different from the definitions in statistical
literature, see, e.g., \cite{CaMoOw,LiOw,Ow,WaKTF}, where superposition 
and truncation dimensions are defined
based on ANOVA decompositions and without any relation to the
integration problem or the error demand $\e$. Moreover, the
dimensions from \cite{CaMoOw, LiOw, WaKTF} depend on specific functions, 
whereas the dimensions in
\cite{KPW15} and in this paper are defined in the worst case sense,
i.e., are relevant to all functions from the space $\calF_\bsgamma$. 

Although the algorithms for constructing $\setU^{\rm q-opt}$ and
$\setU^{\rm opt}$ work for rather general problems and spaces, we
applied them to the integration problem and for
weighted spaces of functions with mixed first order partial
derivatives bounded in $L_p$ norms for $p\in[1,\infty]$.
Such spaces have often been considered (mostly for $p=2$)
when dealing with {\em quasi-Monte Carlo methods}. 

The results depend on how fast the weights converge to zero.
In the experiments, we considered
\[
\gamma_\setu\,=\,\prod_{j\in\setu} j^{-a}\quad\mbox{for\ }a\,=\,2,3,4. 
\]

For $p=1$, the construction from \cite{PlaWas10} is optimal, and it yields
the following results:
\[
d^{\rm sup}(10^{-1})\,=\,\left\{\begin{array}{ll} 1&\mbox{for\ }a=4,\\
2&\mbox{for\ }a=3,\\ 2&\mbox{for\ }a=2, \end{array}\right.
\quad\mbox{and}\quad
|\setU^{\rm PW}(10^{-1})|\,=\,\left\{\begin{array}{ll} 
2 &\mbox{for\ }a=4,\\
4&\mbox{for\ }a=3,\\6 &\mbox{for\ }a=2.\end{array}\right.
\]
\[
d^{\rm sup}(10^{-2})\,=\,\left\{\begin{array}{ll} 2&\mbox{for\ }a=4,\\
2&\mbox{for\ }a=3,\\ 3 &\mbox{for\ }a=2\end{array}\right.
\quad\mbox{and}\quad 
|\setU^{\rm PW}(10^{-2})|\,=\,\left\{\begin{array}{ll} 
6&\mbox{for\ }a=4,\\
8&\mbox{for\ }a=3,\\
22 &\mbox{for\ }a=2.\end{array}\right.
\]
\[
d^{\rm sup}(10^{-3})\,=\,\left\{\begin{array}{ll} 2&\mbox{for\ }a=4,\\
3&\mbox{for\ }a=3,\\  4 &\mbox{for\ }a=2\end{array}\right.
\quad\mbox{and}\quad 
|\setU^{\rm PW}(10^{-3})|\,=\,\left\{\begin{array}{ll} 10
&\mbox{for\ }a=4,\\
22&\mbox{for\ }a=3,\\ 114 &\mbox{for\ }a=2.\end{array}\right.
\]

For $p>1$, $\setU^{\rm PW}$ are no longer optimal; however they are not much
worse than optimal sets when $p$ is relatively close to 1. Moreover,
for all the tests we have performed $d(\setU^{\rm PW}(\e))$ is very
close to the superposition dimension. However the sizes (i.e.,
cardinalities) of $\setU^{\rm opt}(\e)$ and $\setU^{\rm PW}(\e)$ could
be very different, especially for $p=\infty$ and/or small $a$.

For instance, for $p=2$ we have the following results. 
In the case of $\e = 10^{-1}$
\[
d^{\rm sup}(10^{-1})\,\le\,\left\{\begin{array}{ll}
1&\mbox{for\ }a=4,\\
1&\mbox{for\ }a=3,\\
2&\mbox{for\ }a=2,  \end{array}\right.
\]
and
\[
|\setU^{\rm opt}(10^{-1})|\,=\,\left\{\begin{array}{ll}
2&\mbox{for\ }a=4,\\
2 &\mbox{for\ }a=3,\\
4 &\mbox{for\ }a=2,
\end{array}\right.
\quad\mbox{whereas}\quad
|\setU^{\rm PW}(10^{-1})|\,=\,\left\{\begin{array}{ll}
3&\mbox{for\ }a=4,\\
5 &\mbox{for\ }a=3,\\
15 &\mbox{for\ }a=2.
\end{array}\right.
\]

For $\e = 10^{-2}$
\[
d^{\rm sup}(10^{-2})\,\le\,\left\{\begin{array}{ll}
2&\mbox{for\ }a=4,\\
2&\mbox{for\ }a=3,\\
3&\mbox{for\ }a=2  \end{array}\right.
\]
and
\[
|\setU^{\rm opt}(10^{-2})|\,=\,\left\{\begin{array}{ll}
4&\mbox{for\ }a=4,\\
7 &\mbox{for\ }a=3,\\
30 &\mbox{for\ }a=2,
\end{array}\right.
\quad\mbox{whereas}\quad
|\setU^{\rm PW}(10^{-2})|\,=\,\left\{\begin{array}{ll}
8&\mbox{for\ }a=4,\\
18 &\mbox{for\ }a=3,\\
158 &\mbox{for\ }a=2.
\end{array}\right.
\]
And finally, for $\e = 10^{-3}$
\[
  d^{\rm sup}(10^{-3})\,\le\,\left\{\begin{array}{ll} 2&\mbox{for\ }a=4,\\
  3&\mbox{for\ }a=3, \\
  4&\mbox{for\ }a = 2,
  \end{array}\right.
\]
and
\[
|\setU^{\rm opt}(10^{-3})|\,=\,\left\{\begin{array}{ll}
9&\mbox{for\ }a=4,\\ 
24 &\mbox{for\ }a=3,\\
  255&\mbox{for\ }a = 2,
\end{array}\right.
\quad\mbox{whereas}\quad
|\setU^{\rm PW}(10^{-3})|\,=\,\left\{\begin{array}{ll}
20&\mbox{for\ }a=4,\\
70  &\mbox{for\ }a=3,\\
1481&\mbox{for\ }a = 2.
\end{array}\right.
\]
The results for $a=4$ suggest that to achieve an error smaller than $10^{-3}$
it is enough to approximate
$f(\bsx)=\sum_{\setu\subset\bbN_+,|\setu|<\infty}f_\setu(\bsx)$ by 
\[
  f_\emptyset+f_{\{1\}}(x_1)+\cdots+f_{\{5\}}(x_5)+
  f_{\{1,2\}}(x_1,x_2)+f_{\{1,3\}}(x_1,x_3)+f_{\{1,4\}}(x_1,x_4).
\]
As for the quasi-optimal sets, they are the same for $a=4$ and slightly
larger for $a=3,2$: 
\[
|\setU^{\rm q-opt}(10^{-1})|\,=\, 6 \mbox{\ for\ }a=2, 
\quad
|\setU^{\rm q-opt}(10^{-2})|\,=\, 32 \mbox{\ for\ }a=2,
\]
\[
\mbox{and}\quad|
\setU^{\rm q-opt}(10^{-3})|\,=\,
\left\{
\begin{array}{ll}
26 &\mbox{\ for\ } a=3,\\
261 &
\mbox{\ for\ } a = 2.
\end{array}\right.
\]

For $p=\infty$ we have
\[
d^{\rm sup}(10^{-1})\,\leq\,\left\{\begin{array}{ll}
1&\mbox{for\ }a=4,\\
1&\mbox{for\ }a=3,\\
3&\mbox{for\ }a=2,
\end{array}\right.
\]
and
\[
|\setU^{\rm opt}(10^{-1})|\,=\,
\left\{\begin{array}{ll}
2&\mbox{for\ }a=4,\\
3&\mbox{for\ }a=3,\\
33&\mbox{for\ }a=2,
\end{array}\right.
\quad\mbox{whereas}\quad
|\setU^{\rm PW}(10^{-1})|\,=\,
\left\{\begin{array}{ll}
7&\mbox{for\ }a=4,\\
21&\mbox{for\ }a=3,\\
2358&\mbox{for\ }a=2.
\end{array}\right.
\]
Now for $\e = 10^{-2}$
\[
d^{\rm sup}(10^{-2})\,\leq\,\left\{\begin{array}{ll}
2&\mbox{for\ }a=4,\\ 
2&\mbox{for\ }a=3,\\
4&\mbox{for\ }a=2,
\end{array}\right.
\]
and
\[
|\setU^{\rm opt}(10^{-2})|\,=\,\left\{\begin{array}{ll}
5&\mbox{for\ }a=4,\\ 
15&\mbox{for\ }a=3,\\
1346&\mbox{for\ }a=2,
\end{array}\right.
\quad\mbox{whereas}\quad
|\setU^{\rm PW}(10^{-2})|\,=\,\left\{\begin{array}{ll}
21&\mbox{for\ }a=4,\\ 149 &\mbox{for\ }a=3\\
120,935&\mbox{for\ }a=2.
\end{array}\right.
\]
For $\e = 10^{-3}$
\[
d^{\rm sup}(10^{-3})\,\leq\,\left\{\begin{array}{ll}
2&\mbox{for\ }a=4,\\ 
3&\mbox{for\ }a=3,\\
6&\mbox{for\ }a=2,
\end{array}\right.
\]
and
\[
|\setU^{\rm opt}(10^{-3})|\,=\,\left\{\begin{array}{ll}
15&\mbox{for\ }a=4,\\
83&\mbox{for\ }a=3,\\
45,446&\mbox{for\ }a=2,
\end{array}\right.\quad\mbox{whereas}\quad
|\setU^{\rm PW}(10^{-3})|\,=\,\left\{\begin{array}{ll}
72&\mbox{for\ }a=4,\\
923 &\mbox{for\ }a=3.
\end{array}\right.
\]

For the tests above, the quasi-optimal active set was different from the
corresponding optimal active set in the following cases only:
\[
|\setU^{\rm q-opt}(10^{-1})|\,=\, 38 \mbox{\ for\ }a=2, \quad
|\setU^{\rm q-opt}(10^{-2})|\,=\, 1904 \mbox{\ for\ }a=2, 
\]
\[
\mbox{and}\quad|
\setU^{\rm q-opt}(10^{-3})|\,=\,
\left\{
\begin{array}{ll}
92 &\mbox{\ for\ } a=3,\\
52,159 &
\mbox{\ for\ } a = 2.
\end{array}\right.
\]

A collection of the active sets constructed above have been listed in full in the
Appendix.

Our algorithms can also be used to construct the active sets 
where, instead of the standard worst case error, the normalized 
worst case error is used. More precisely, for the normalized worst
case error
we would like to have sets $\setU_{\rm norm}(\e)$ such that the integral of 
$\sum_{\setu\notin\setU_{\rm norm}(\e)} f_\setu$ is bounded by 
\[
  \e\,\|\calS\|\,\|f\|\quad\mbox{for all\ }f\in\calF_\bsgamma,
\]
where $\|\calS\|$ is the norm of the integration operator. Since in
our case $\|\calS\|\ge 1$, the corresponding active sets 
$\setU_{\rm norm}^{\rm X}(\e)$ are subsets of $\setU^{\rm X}(\e)$ 
for ${\rm X}\in\{{\rm opt,\,q-opt,\,PW}\}$ and could be even smaller.

\section{Basic Definitions}
We provide in this section basic concepts and definitions for special
spaces of functions that are very often assumed in the 
literature, especially in the context of quasi-Monte Carlo methods.
The presented algorithms can easily be modified to more
general spaces. 

\subsection{$\bsgamma$-Weighted Spaces}
We follow here \cite{HefRitWas15}.
For $D=[0,1]$,  let $\calD=D^{\bbN_+}$ be the set of sequences (points)
$\bsx=[x_1,x_2,\dots]$ with $x_i\in D$. Here $\bbN_+$ is the set of 
positive integers and we will use $\setu$ and $\setv$ to denote finite
subsets of $\bbN_+$. We will also use the following notation:
For $\bsx\in\calD$ and $\setu$,
by $[\bsx_\setu;\bszero_{\setu^c}]$ we denote the point in $\calD$ such that
\[
[\bsx_\setu;\bszero_{\setu^c}]\,=\,[y_1,y_2,\dots]\quad\mbox{with}\quad
y_j\,=\,\left\{\begin{array}{ll} x_j&\mbox{if\ }j\in\setu,\\
0& \mbox{if\ }j\notin\setu.
\end{array}\right.
\]
Next,
\[
f^{(\setu)}\,=\,\prod_{j\in\setu}\frac{\partial}{\partial x_j} f.
\]

For given
\[
p\in[1,\infty],
\]
let $\calF_{\bsgamma,p}$ be the Banach space of functions defined on
$\calD$ with the following norm
\[ 
\|f\|_{\calF_{\bsgamma,p}}\,=\,
\left(\sum_{\setu\subset\bbN_+,|\setu|<\infty} \gamma_\setu^{-p}\,
\|f^{(\setu)}([\cdot_\setu;\bszero_{\setu^c}])\|_{L_p}^p\right)^{1/p}.
\]
Of course, for $p=\infty$,
\[\|f\|_{\calF_{\bsgamma,p}}\,=\,
\sup_{\setu\subset\bbN_+,|\setu|<\infty}\frac{\|f^{(\setu)}([\cdot_\setu;
    \bszero_{\setu^c}])\|_{L_\infty}}{\gamma_\setu}.
\]
We assume that the numbers $\gamma_\setu$ (called {\em weights}) are of 
product form (see \cite{SloWoz98})
\begin{equation}\label{prod-w}
\gamma_\setu\,=\,\prod_{j\in\setu}\frac{c}{j^a}\quad\mbox{for
  positive $a$ and $c$}.
\end{equation}
In general choosing the weights (in our case choosing $a$ and $c$) for
a specific integral or application is a difficult problem which we do not 
attempt to address here.
We assume that the parameters $a$ and $c$ are given with the problem.

It was shown in \cite{HefRitWas15} that any $f\in\calF_{\bsgamma,p}$
admits a unique decomposition, called the {\em anchored decomposition},
\[f\,=\,\sum_{\setu\subset\bbN_+,|\setu|<\infty}f_\setu
\]

with $f_\setu$ given by 
\[
f_\setu(\bsx)\,=\,T_\setu(h_\setu)(\bsx)
\,:=\,\int_{D^{|\setu|}}h_\setu(\bst)\,
\prod_{j\in\setu}(x_j-t_j)^0_+\rd\bst\quad\mbox{for some $h_\setu\in 
L_p(D^{|\setu|})$},
\]
where $(x_j-t_j)^0_+$ is $1$ if $x_j>t_j$ and $0$ otherwise. 
The functions $f_\setu$ belong to the following Banach spaces 
$F_\setu$ 
\[F_\setu\,=\,T_\setu(L_p)\quad\mbox{and}\quad\|f_\setu\|_{F_\setu}
\,=\,\|f_\setu^{(\setu)}\|_{L_p}.
\]
Of course, $F_{\emptyset}$ is the space of constant functions
with the absolute value as its norm. 
The spaces $F_\setu$ for $\setu\not=\emptyset$ are anchored at $0$ since
$f_\setu(\bsx)=0$ if there is $j\in\setu$ with $x_j=0$. This is why
\[
f^{(\setu)}([\cdot_\setu;\bszero_{\setu^c}])\,=\,f_\setu^{(\setu)}
\qquad\mbox{and}\qquad 
\|f\|_{\calF_{\bsgamma,p}}\,=\,
\left(\sum_{\setu\subset\bbN_+,|\setu|<\infty}
\gamma_\setu^{-p}\,\|f_\setu\|_{F_\setu}^p\right)^{1/p}.
\]

The space $\calF_{\bsgamma,p}$ contains in particular the following
class of functions.

\begin{example}
For a smooth function $g:\bbR\to\bbR$ 
and fast decaying numbers $a_1,a_2,\dots$, consider 
\begin{equation}\label{exmp}
   f(\bsx)\,=\,g\left(\sum_{j=1}^\infty x_j\,
     a_j\right)\quad\mbox{for\ }x_j\in D.
\end{equation}
Then 
\[
   f^{(\setu)}([\bsx_\setu;\bszero_{\setu^c}])\,=\,g^{(|\setu|)}
   \left(\sum_{j\in\setu}x_j\,
     a_j\right)\, \prod_{j\in\setu}a_j. 
\]
Hence $f\in\calF_{\bsgamma,p}$ if 
the derivatives of $g$ and the coefficients $a_j$ satisfy
\[
  \left(\sum_{\setu\subset\bbN_+,|\setu|<\infty}
   \frac{\prod_{j\in\setu}|a_j|^p}{\gamma_\setu^p}\,
     \int_{D^{|\setu|}}\left|g^{(|\setu|)}\left(\sum_{j\in\setu}
   x_j\, a_j\right)\right|^p\rd \bsx_\setu\right)^{1/p}\,<\,\infty.
\]
\end{example}

\subsection{Integration Problem}\label{sec:int}
Consider the following integration functional
\[
\calS:\calF_{\bsgamma, p}\to\bbR
\]
given by
\[
\calS(f)\,=\,\lim_{s\to\infty}\int_{D^s}f(x_1,\dots,x_s,0,\dots,0)
\rd[x_1,\dots,x_s].
\]
Let $p^*$ denote the conjugate of $p$,
\[
  \frac1p+\frac1{p^*}\,=\,1.
\]
We assume that
\begin{equation}\label{ass-2}
  \left(\sum_{\setu\subset\bbN_+,|\setu|<\infty}
  \gamma_\setu^{p^*}\,(p^*+1)^{-|\setu|}\right)^{1/p^*}
  \,<\,\infty
\end{equation}
since the left hand side of \eqref{ass-2} is the
norm of $\calS$, i.e., \eqref{ass-2} is a necessary and sufficient
condition for continuity of $\calS$. Indeed, letting $S_\setu$ be the
restriction of $\calS$ to $F_\setu$, we have that
\[
\|S_\setu\|_{F_\setu}\,=\,\sup_{\|f_\setu\|_{F_\setu}=1} S_\setu(f_\setu)
\,=\,\frac1{(p^*+1)^{|\setu|\,/p^*}}
\]
which, with an application of H\"older's inequality, yields \eqref{ass-2}.
For product weights of the form \eqref{prod-w}, we have
\[
\|\calS\|\,=\,\left(\sum_\setu \gamma_\setu^{p^*}\,(p^*+1)^{-|\setu|}
\right)^{1/p^*}\,=\,\prod_{j=1}^\infty\left(1+\frac{c^{p^*}}{j^{a\,p^*}\,
  (p^*+1)}\right)^{1/p^*}.
\]
Hence for product weights, \eqref{ass-2} is equivalent to 
$a\,>\,1/{p^*}$. For the remainder of the paper it is assumed
that $a \,>\, 1/p^*$.

A very important part of the {\em Multivariate Decomposition Method}
({\em MDM} for short) is a construction of {\em active sets}
$\setU(\e)$, i.e., sets that satisfy 
\begin{equation}\label{mod-red}
  \left|\calS\left(\sum_{\setu\notin\setU(\e)}f_\setu\right)\right|\,
  \le\,\e\,
  \left\|\sum_{\setu\notin\setU(\e)}f_\setu\right\|_{\calF_{\bsgamma,p}}
\quad\mbox{for all\ } f\in\calF_{\bsgamma,p}.
\end{equation}
The essence of \eqref{mod-red} is that, when approximating $\calS(f)$, 
it is enough to restrict the attention to functions 
\[
  \sum_{\setu\in\setU(\e)} f_\setu,
\]
since any algorithm approximating
$\sum_{\setu\in\setU(\e)}S_\setu(f_\setu)$
with the worst case error on $\bigoplus_{\setu\in\setU(\e)}F_\setu$
bounded by $\e$ has its worst case
error on the whole space $\calF_{\bsgamma,p}$ bounded by
\[
    2^{1/p^*}\,\e. 
\]
The factor of $2^{1/p^*}$ is the result of applying H\"older's inequality,
see, e.g., \cite{KPW15}.
Clearly, there are many sets satisfying \eqref{mod-red}, and we would like
to construct possibly small active sets.
\begin{definition}
We say that an active set, denoted by $\setU^{\rm opt}(\e)$ is
{\em optimal}, if
\[
|\setU^{\rm opt}(\e)|\,=\,\min\{|\setU(\e)|\ :\ \setU(\e)\mbox{\
    satisfies\ }\eqref{mod-red}\}.
\]
We also define the $\e$-{\em superposition dimension} as the
smallest $d(\setU(\e))$ among all active sets,
\[
  d^{\rm sup}(\e)\,:=\,\min\left\{d(\setU(\e))\ :\ \setU(\e)\ 
\mbox{satisfies \eqref{mod-red}}\right\}. 
\]
\end{definition}

\section{Constructing Active Sets $\setU(\e)$}
A construction of active sets was first proposed in \cite{PlaWas10}. 
The corresponding sets will be denoted by $\setU^{\rm PW}(\e)$. 
It was shown there that 
\[
  d(\setU^{\rm PW}(\e))\,=\,O\left(\frac{\ln(1/\e)}{\ln(\ln(1/\e))}\right)
\quad\mbox{as\ }\e\,\to\,0 
\]
and that the cardinality of $\setU^{\rm PW}(\e)$ is polynomial in $1/\e$. 
It is easy to see that for $p=1$, $\setU^{\rm PW}(\e)$ are optimal. 
However, as we shall see, their size might be too large
for large values of $p$, especially for  $p=\infty$. 

Let $\setU$ be a set of subsets $\setu$. Then 
\begin{eqnarray*}
  \left|\calS\left(\sum_{\setu\notin\setU}f_\setu\right)\right|
  &\le&\sum_{\setu\notin\setU}\|f_\setu\|_{F_\setu}\,
    \|S_\setu\|_{F_\setu}\,=\,
  \sum_{\setu\notin\setU}\frac{\|f_\setu\|_{F_\setu}}{\gamma_\setu}\,
   \gamma_\setu\,\|S_\setu\|_{F_\setu}\\
  &\le&\left\|\sum_{\setu\notin\setU}f_\setu\right\|_{\calF_{\bsgamma,p}}\,
   \left(\sum_{\setu\notin\setU}\gamma_\setu^{p^*}\,
  \|S_\setu\|_{F_\setu}^{p^*}\right)^{1/p^*}. 
\end{eqnarray*}
Hence we are looking for $\setU(\e)$ such that 
\begin{equation}\label{mod-red2}
   \left(\sum_{\setu\notin\setU(\e)}\gamma_\setu^{p^*}\,
  \|S_\setu\|_{F_\setu}^{p^*}\right)^{1/p^*}\,\le\,\e.
\end{equation}
Since H\"older's inequality is sharp, \eqref{mod-red2} is equivalent
to \eqref{mod-red}.

For the sake of completeness, we recall the construction for 
$p=1$, see \cite{PlaWas10}. 

\subsection{Case $p=1$}
For $p=1$, we have $p^*=\infty$ and $\|S_\setu\|_{F_\setu}=1$ for 
all $\setu$. Hence
\[
 \left(\sum_{\setu\notin\setU(\e)}\gamma_\setu^{p^*}\,
\|S_\setu\|^{p^*}_{F_\setu}\right)^{1/p^*}\,=\,\sup_{\setu\notin\setU(\e)}
 \gamma_\setu
\]
which for product weights reduces to
$\sup_{\setu\notin\setU(\e)}\prod_{j\in\setu}c\,j^{-a}$. 
Therefore 
\begin{equation}\label{setu-inf}
  \setU^{\rm PW}(\e)
  \,=\,\left\{\setu\ :\ \prod_{j\in\setu}\frac{c}{j^a}\,>\,\e
  \right\}.
\end{equation}
It is easy to see that $\setU^{\rm PW}(\e)$ is 
the smallest set satisfying  \eqref{mod-red}, i.e., it is a subset
of any $\setU(\e)$ satisfying \eqref{mod-red2}. 

The examples of $\setU^{\rm PW}$ for specific values
of $a$ and $\e$ are presented in the Appendix. 
For simplicity we use $c=1$ there. 

\subsection{Case of $p>1$}
For $p>1$, the conjugate $p^*$ is finite and the construction of
$\setU(\e)$ is more complicated.

We begin by recalling the construction of $\setU^{\rm PW}(\e)$ in 
\cite{PlaWas10}. 
To simplify the notation, let
\[
\overline{\gamma}_\setu\,=\,\frac{\gamma_\setu^{p^*}}{(p^*+1)^{|\setu|}}
\,=\,\left(\frac{c^{p^*}}{p^*+1}\right)^{|\setu|}\,
\prod_{j\in\setu} j^{-a\,p^*}. 
\]
For given $\e$ and $p$, a special threshold is
computed and all $\setu$ with $\overline{\gamma}_\setu$ 
exceeding the threshold are included in the active set. 
More precisely, for $t\in(1/(ap^*),1)$ a threshold is given by
\[
{\rm Threshold}(\e,t)\,=\,\left(
\frac{\e^{p^*}}{\sum_{\setu\subset\bbN_+,|\setu|<\infty}
  \overline{\gamma}_\setu^{\,t}}\right)^{1/(1-t)}.
\]
Note that the interval $(1/(ap^*), 1)$ is non-empty by the assumption
that $a \,>\,1/p^*$ introduced in Section~\ref{sec:int}.
In our numerical experiments we approximated the sum of
$\overline{\gamma}_\setu^{\,t}$ for $t=i/40$ ($39\ge i>40/(ap^*)$)
and selected the value which resulted in the largest ${\rm Threshold}(\e, t)$. 
The approximations are calculated in a similar way as the computation of
$A_s$ explained later (see \eqref{eq:A_s}).

Clearly, $\setU^{\rm PW}(\e)$ contains a number
of $\setu$'s with the largest $\overline{\gamma}_\setu$; 
however, the number of them could be much larger than needed. 

The collection of those $\setu$ with the largest 
$\overline{\gamma}_\setu$ that are necessary 
for \eqref{mod-red} would result in the optimal set $\setU^{\rm  opt}(\e)$.
Since this optimal set is always a subset of $\setU^{\rm  PW}(\e)$, 
the good property
\[
d(\setU^{\rm opt}(\e))\,=\,O\left(\frac{\ln(1/\e)}{\ln(\ln(1/\e))}\right)
\quad\mbox{as\ }\e\to0,
\] 
is preserved. 
More precisely, let $(\setu_j)_{j\in\bbN_+}$ be a sequence of
all subsets $\setu$ ordered so that
\[
\overline{\gamma}_{\setu_j}\,\ge\,\overline{\gamma}_{\setu_{j+1}}
\quad j=1,2,\dots.
\]
Then
\[
\setU^{\rm opt}(\e)\,=\,\left\{\setu_1,\dots,\setu_k\right\}
\]
with $k=k(\e)$ such that  
\[
\|\calS\|^{p^*}-\sum_{j=1}^k \overline{\gamma}_{\setu_j}
\,\le\,\e^{p^*}\,<\,
\|\calS\|^{p^*}-\sum_{j=1}^{k-1} \overline{\gamma}_{\setu_j}.
\]
The problem with this approach is that we do not know {\em a priori}
the number $k=k(\e)$ and ordering a large number of
$\overline{\gamma}_\setu$ might be too expensive. This is why the
numbers $\overline{\gamma}_\setu$ will be ordered on-line. 
Actually, we propose two ways of constructing active sets. The first
and simpler one produces what we call, {\em quasi-optimal} sets
$\setU^{\rm q-opt}(\e)$ and it uses a partial ordering of
$\overline{\gamma}_\setu$. The second
one, uses ordering of $\overline{\gamma}_\setu$ and produces optimal
sets $\setU^{\rm opt}(\e)$. 
However, as we will see the difference between both sets is very
small; sometimes these sets are equal. 

The numbers $\overline{\gamma}_\setu$ have the following properties
that are crucial for our construction of quasi-optimal and optimal 
sets. 
Let $\ell$ be a given cardinality. In what follows we will write
\[
\setu\,=\,\left\{u_1,\dots,u_\ell\right\},\quad\mbox{where}\quad
u_1\,<\,\cdots\,<\,u_\ell.
\]
The first property is: If 
\[
\setu\,=\,\{u_1,\dots,u_\ell\}\mbox{\ and\ }
\setv\,=\,\{v_1,\dots,v_\ell\}
\mbox{\ with\ }v_j\,\ge\,u_j\mbox{ for all\ }j
\]
then
\[
\overline{\gamma}_\setu\,\ge\,\overline{\gamma}_\setv.
\]
The other property is: For $\ell+1\ge c^{1/a}$,
\[
\overline{\gamma}_{\{u_1,\dots,u_\ell\}} \,\ge\,
\overline{\gamma}_{\{u_1,\dots,u_\ell,u_{\ell+1}\}}.
\]

We are ready to describe the constructions of active sets. First we need
to approximate 
\[
A\,=\,\sum_{\setu\subset\bbN_+,|\setu|<\infty}\overline{\gamma}_\setu
\]
from above and with the relative error significantly smaller than
$\e^{p^*}$. 
This can be done as follows. For a large natural number $s$
\begin{eqnarray}
  A&=&\exp\left(\ln\left(\prod_{j=s+1}^\infty\left(1+\frac{(c/j^a)^{p^*}}
  {p^*+1}\right)\right)\right)\,
\prod_{j=1}^s\left(1+\frac{(c/j^a)^{p^*}}{p^*+1}\right)\nonumber\\
&\le& \exp\left(\frac{c^{p^*}}{p^*+1}\sum_{j=s+1}^\infty
  j^{-a\,p^*}\right)\,
\prod_{j=1}^s\left(1+\frac{(c/j^a)^{p^*}}{p^*+1}\right)\nonumber\\
&\le& \exp\left(\frac{c^{p^*}}{p^*+1}\int_{s+1/2}^\infty
  x^{-a\,p^*}\rd x\right)\,
\prod_{j=1}^s\left(1+\frac{(c/j^a)^{p^*}}{p^*+1}\right)\nonumber\\
&=&\exp\left(\frac{c^{p^*}}{(p^*+1)\,(a\,p^*-1)\,(s+1/2)^{a\,p^*-1}}\right)
\,\prod_{j=1}^s\left(1+\frac{(c/j^a)^{p^*}}{p^*+1}\right)\,=:\,A_s.
\label{eq:A_s}
\end{eqnarray}
It is easy to see that the relative error between $A$ and its
approximation $A_s$ is proportional to
$1/s^{2a\,p^*-2}$ with the asymptotic constant $c^{p^*}/((p^*+1)\,
2^{ap^*-1})\,\prod_{j=1}^\infty(1+(c/j^a)^{p^*}/(p^*+1))$. 

A general idea of our construction is to
select sets $\setu$ with large $\overline{\gamma}_\setu$ and subtract
$\overline{\gamma}_\setu$ from $A_s$. This is repeated until
$A_s$ is reduced to or below $\e^{p^*}$. 

More specifically, 
consider a partition of $\bbR_+$ into intervals $I_i$ such that the
numbers in $I_j$ are greater than those in $I_{j+1}$. For simplicity, we
used 
\[
I_1\,=\,[10^{-1},\infty),\quad\mbox{and}\quad
  I_j\,=\,[10^{-j},10^{-j + 1})\quad\mbox{for\ }j=2,3,\dots.
\]
in our numerical experiments when constructing quasi-optimal sets. 
However, we think that a better partition is
possible, especially when constructing optimal active sets.
We also associate with every interval a list $L_j$ that
contains those $\setu$ for which $\overline{\gamma}_\setu$
has been subtracted from $A_s$ in $j$th step.

In the first $j=1$ step, add the empty set to $L_1$ and subtract
$\overline{\gamma}_\emptyset=1$
from $A_s$. If the new $A_s$ satisfies $A_s \le\e^{p^*}$, then terminate.
Otherwise consider non-empty sets $\setu$ in the order of increasing
cardinalities. Hence start with singleton sets $\setu=\{i\}$ for
$i=1,\dots,k$, 
where $k$ is the largest integer such that $\overline{\gamma}_{\{k\}}$
is in $I_1$. Place $\{k+1\}$ into list $L_2$, and start subtracting
from $A_s$ the values $\overline{\gamma}_{\{i\}}$ and store $\{i\}$
in $L_1$ until either the
difference becomes less than or equal to $\e^{p^*}$, in which case we terminate, or $i=k$.
Next repeat the same for sets of cardinality $2$, starting with
sets $\{1,i\}$ for $i\le k$, where now $k$ is the largest integer such that
$\overline{\gamma}_{\{1,k\}}\in I_1$. Store $\{1,k+1\}$ in $L_2$.
Next consider sets $\{2,i\}$, $\{3,i\}$, etc., until either all cardinality 2 
sets corresponding to the current interval have been visited or the new value of
$A_s$ is $\le\e^{p^*}$, in which case we terminate. Continue working through
the sets in order of increasing cardinality $\ell$ 
until $\overline{\gamma}_{\{1,2,\dots,\ell\}}\notin I_1$ and $\ell\ge c$
(of course, $\ell$ is always at least $c$ for $c\le1$). 
Then move to step $j=2$.
The procedure in this (and later) steps is very similar except that
for a fixed cardinality of $\setu$, we check if any such set
has already been placed in $L_2$ in the 1st step. If it has, we start
working with such sets first. For instance, for cardinality
$1$, if $\{k+1\}\in L_2$, then we begin with sets $\{i\}$ for $i\ge k+1$.
For cardinality $2$, if $\{i_i,i_2\}$ has been placed in $L_2$, then we
inspect sets $\{i_1,i\}$ for $i\ge i_2$, before any other sets of
cardinality two are considered.
Once we find $\overline{\gamma}_{\{i,i+1\}}\notin I_2$,
we store $\{i,i+1\}$ in $L_3$ and proceed to sets of cardinality 3, etc. 

At the very end, $\setU^{\rm q-opt}(\e)$ consists of all subsets
$\setu$ whose values $\overline{\gamma}_\setu$ were subtracted from $A_s$. 
The main procedure is outlined in Algorithm~\ref{alg:q-opt}.
In all of the algorithms $j_{\max}$ and $\ell_{\max}$ are computational thresholds
denoting, respectively, the maximum number of intervals to be searched through and 
the maximum allowed cardinality of sets.

To search through the sets in a systematic way, we must keep track of the 
current set, $\setu$, and the index, $i$, that we are incrementing from.
The subroutine \textbf{increment-$\setu$} outlined below in Algorithm~\ref{alg:incr-u}
details how to increment $\setu$ from index $i$.

\begin{algorithm}[!h]
\caption{(Subroutine: \textbf{increment-$\setu$})}\label{alg:incr-u}
\textbf{inputs}: $\setu$, $i$\\
\textbf{output}: $\setu$
\begin{algorithmic}[1]
\State $u_i \leftarrow u_i + 1$
\Comment updating $u_i$ first
\For{{$r = i + 1,i + 2, \ldots, |\setu|$}}
\Comment incrementing $\setu$ from index $i + 1$
\State $u_r \leftarrow u_i + r - i$
\EndFor
\State return $\setu$
\end{algorithmic}
\end{algorithm}

\begin{algorithm}[!h]
\caption{(Constructing the quasi-optimal active set)}\label{alg:q-opt}
\textbf{inputs}: $\e$, $p^*$, $s$, $(\bar{\gamma}_\setu)_{|\setu| < \infty}$, $(I_j)_{j = 1}^{j_{\max}}$\\
\textbf{output}: $\setU^\mathrm{q-opt}(\e)$
\begin{algorithmic}[1]
\State $L_j \leftarrow \emptyset$ for all $j = 1, 2, \ldots, j_{\max}$
\Comment initialising
\State $\setU^{\rm q-opt}(\e) \leftarrow \{\emptyset\}$
\State $T \leftarrow A_s - \e^{p^*} - \bar{\gamma}_\emptyset$
\Comment $T$ tracks difference between $A_s - \e^{p^*}$ and weights
\State \textbf{if} $T \leq 0$ \textbf{then} return $\setU^\mathrm{q-opt}(\e)$
\Comment quasi-optimal active set is complete
\For{{$j = 1, 2, \ldots, j_{\max}$}}
\Comment looping over intervals
\State \Comment first handle sets found at previous step
\State $(\setU^\mathrm{q-opt}(\e), T, \ell_\mathrm{next}, L_{j + 1})
\leftarrow $\textbf{q-opt-search}($\setU^\mathrm{q-opt}(\e), T, 
(\bar{\gamma}_\setu)_{|\setu| < \infty}, L_j, L_{j + 1}, I_j)$
\State \textbf{if} $T \leq 0$ \textbf{then} return $\setU^\mathrm{q-opt}(\e)$
\Comment quasi-optimal active set is complete
\For{{$\ell = \ell_\mathrm{next}, \ell_\mathrm{next} + 1, \ldots, \ell_{\max}$}}
\Comment search through unvisited sets
\State $\setu = \{1, 2, \ldots, \ell\}$
\State $i \leftarrow \ell$
\Comment $i$ keeps track of index to increment $\setu$ from
\State \textbf{if} $\bar{\gamma}_\setu \notin I_j$ \textbf{and} $\ell \geq c$
\textbf{then} break
\Comment no more $\setu$ with $\bar{\gamma}_{\setu} \in I_j$
\While{{$i > 0$}}
\Comment when $i = 0$ there are no more $\setu$ of cardinality $\ell$
\If{{$\bar{\gamma}_\setu \in I_j$}}
\State add $\setu$ to $\setU^\mathrm{q-opt}(\e)$
\State $T \leftarrow T - \bar{\gamma}_\setu$
\State \textbf{if} $T \leq 0$ \textbf{then} return $\setU^\mathrm{q-opt}(\e)$
\Comment quasi-optimal set is complete
\State $i \leftarrow \ell$
\Comment continue incrementing from last index
\Else
\State add $\setu$ to $L_{j + 1}$
\State $i \leftarrow i - 1$
\Comment start incrementing from lower index
\State \textbf{if} $i = 0$ \textbf{then} break
\Comment go to next cardinality
\EndIf
\State $\setu \leftarrow $\textbf{increment-$\setu$}$(\setu, i)$
\EndWhile
\EndFor
\EndFor
\end{algorithmic}
\end{algorithm}
\pagebreak
To make the presentation clearer Algorithm~\ref{alg:q-opt} is broken into two parts: 
First, we search starting from the sets found in the previous interval, which is handled
by the subroutine \textbf{q-opt-search} in Algorithm~\ref{alg:search-q}.
Then we continue searching through sets in order of increasing cardinality (line 9)
starting where \textbf{q-opt-search} finished, at cardinality $\ell_\mathrm{next}$.
The basic search structure is the same, however in \textbf{q-opt-search} 
each set we visit is checked to reduce multiple visits to a single set and
ensure that the same set is not added to $\setU^\mathrm{q-opt}(\e)$ more than once.

The notation
\[
(\setU^\mathrm{q-opt}(\e), T, \ell_\mathrm{next}, L_{j + 1})
\,\leftarrow\,
\text{\textbf{q-opt-search}}
(\setU^\mathrm{q-opt}(\e),T, (\bar{\gamma}_\setu)_{|\setu| < \infty}, L_j, L_{j + 1}, I_j)\,,
\]
denotes that we call \textbf{q-opt-search} with inputs 
$\setU^\mathrm{q-opt}(\e)$, $T$,
$(\bar{\gamma}_\setu)_{|\setu| < \infty}$, $L_j$, $L_{j + 1}$, $I_j$ and then 
use the output to update $\setU^\mathrm{q-opt}(\e)$, $T$, $\ell_\mathrm{next}$ and $L_{j + 1}$.

\begin{algorithm}[!h]
\caption{(Subroutine: \textbf{q-opt-search})}\label{alg:search-q}
\textbf{inputs}: $\setU^\mathrm{q-opt}(\e)$, $T$, $(\bar{\gamma}_\setu)_{|\setu| < \infty}$, 
$L_j$, $L_{j + 1}$, $I_j$\\
\textbf{outputs}: $\setU^\mathrm{q-opt}(\e)$, $\ell_\mathrm{next}$, $T$, $L_{j+ 1}$
\begin{algorithmic}[1]
\State $\ell_\mathrm{next} = 1$
\For{{$\setu \in L_j$}}
\State $i \leftarrow |\setu|$
\While{{$i > 0$}}
\State $L_j \leftarrow L_j \setminus \setu$
\Comment reducing the double-handling of sets
\If{{$\bar{\gamma}_\setu \in I_j$}}
\State \textbf{if} $\setu \in \setU^\mathrm{q-opt}(\e)$ \textbf{then} break
\Comment already visited $\setu$ and any future increments
\State add $\setu$ to $\setU^\mathrm{q-opt}(\e)$
\State $T \leftarrow T - \bar{\gamma}_\setu$
\State \textbf{if} $T \leq 0$ \textbf{then} return 
$(\setU^\mathrm{q-opt}(\e), \ell_\mathrm{next}, T, L_{j + 1})$
\State $i = |\setu|$
\Else
\State add $\setu$ to $L_{j + 1}$
\Comment $\setu$ to be checked first in next interval
\State $i \leftarrow i - 1$
\State \textbf{if} $i = 0$ \textbf{then} break
\Comment go to next $\setu \in L_j$
\EndIf
\State $\setu \leftarrow $\textbf{increment-$\setu$}$(\setu, i)$
\EndWhile
\State $\ell_\mathrm{next} = |\setu| + 1$
\Comment main search will start at cardinality $|\setu| + 1$
\EndFor
\State return $(\setU^\mathrm{q-opt}(\e), \ell_\mathrm{next}, T, L_{j + 1})$
\end{algorithmic}
\end{algorithm}

The construction of optimal active sets is very similar. The main 
difference is that in the $j$th step, we first create the list $L_{j}^{\rm unsorted}$,
order its elements $\setu\in L_{j}^{\rm unsorted}$ according to decreasing values of 
$\overline{\gamma}_{\setu}$, and next start subtracting the values 
$\overline{\gamma}_\setu$ from $A_s$.  
Again the lists $L_j$ will hold the sets visited in
the previous interval.

In fact, if we do not care whether
or not all of the sets are ordered but only that $\setU^{\rm opt}(\e)$ consists
of the sets with the largest weights, then 
we only need to sort the sets which
come from the final interval. This is because at the previous intervals all of the
sets will need to be added to $\setU^{\rm opt}(\e)$, regardless of sorting. To do this in
practice, for each interval we store the sum of all the weights corresponding
to that interval. 
In Algorithm~\ref{alg:opt} we denote this by $T_j$. 
At the end of the $j$th
step, we check whether $I_j$ is the final interval, i.e., if $A_s - \sum_{i = 1}^jT_i \leq \varepsilon^{p^\star}$,
if so we sort the sets and add them one-by-one until the active set is complete.
Otherwise we add all of the sets in $L_j^{\rm unsorted}$ to $\setU^{\rm opt}(\e)$  and
go to the next interval. 
For completeness, the construction of optimal active sets is detailed 
separately below in Algorithm~\ref{alg:opt} and the subroutine
\textbf{opt-search} in Algorithm~\ref{alg:search-opt}.

\begin{algorithm}[h!]
\caption{(Constructing the optimal active set)}\label{alg:opt}
\textbf{inputs}: $\e$, $p^*$, $s$, $(\bar{\gamma}_\setu)_{|\setu| < \infty}$, $(I_j)_{j = 1}^{j_{\max}}$\\
\textbf{output}: $\setU^\mathrm{opt}(\e)$
\begin{algorithmic}[1]
\State $T_j \leftarrow 0$, $L_j^\mathrm{unsorted} \leftarrow \emptyset$ and $L_j \leftarrow \emptyset$ for all $j = 1, 2, \ldots, j_{\max}$
\Comment initialising
\State $\setU^{\rm opt}(\e) \leftarrow \{\emptyset\}$
\State $T \leftarrow A_s - \e^{p^*} - \bar{\gamma}_\emptyset$
\Comment $T$ tracks difference between $A_s - \e^{p^*}$ and weights
\State \textbf{if} $T \leq 0$ \textbf{then} return $\setU^\mathrm{opt}(\e)$
\Comment optimal active set is complete
\For{{$j = 1, 2, \ldots, j_{\max}$}}
\Comment looping over intervals
\State \Comment first handle sets found at previous step
\State $(\ell_\mathrm{next}, T_j, L_j^\mathrm{unsorted}, L_{j + 1})
\leftarrow $\textbf{opt-search}($\setU^\mathrm{opt}(\e), T_j, 
(\bar{\gamma}_\setu)_{|\setu| < \infty}, L_j^\mathrm{unsorted}, L_{j + 1}, I_j)$
\For{{$\ell = \ell_\mathrm{next}, \ell_\mathrm{next} + 1, \ldots, \ell_{\max}$}}
\Comment search through unvisited sets
\State $\setu = \{1, 2, \ldots, \ell\}$
\State $i \leftarrow \ell$
\Comment $i$ keeps track of index to increment $\setu$ from
\State \textbf{if} $\bar{\gamma}_\setu \notin I_j$ \textbf{and} $\ell \geq c$
\textbf{then} break
\Comment no more $\setu$ with $\bar{\gamma}_{\setu} \in I_j$
\While{{$i > 0$}}
\Comment when $i = 0$ there are no more $\setu$ of cardinality $\ell$
\If{{$\bar{\gamma}_\setu \in I_j$}}
\State add $\setu$ to $L_j^\mathrm{unsorted}$
\State $T_j \leftarrow T_j + \bar{\gamma}_\setu$
\State $i \leftarrow \ell$
\Comment continue incrementing from last index
\Else
\State add $\setu$ to $L_{j + 1}$
\State $i \leftarrow i - 1$
\Comment start incrementing from lower index
\State \textbf{if} $i = 0$ \textbf{then} break
\Comment go to next cardinality
\EndIf
\State $\setu \leftarrow $\textbf{increment-$\setu$}$(\setu, i)$
\EndWhile
\EndFor
\If{{$T_j \geq T$}}
\Comment sorting step, first check if $I_j$ is the last interval
\State sort $L_j^\mathrm{unsorted}$
\For{{$\setu \in L_j^\mathrm{sorted}$}}
\Comment add sorted sets until active set is complete
\State add $\setu$ to $\setU^\mathrm{opt}(\e)$
\State $T \leftarrow T - \bar{\gamma}_\setu$
\State \textbf{if} $T \leq 0$ \textbf{then} return $\setU^\mathrm{opt}(\e)$
\Comment optimal active set is complete
\EndFor
\Else 
\Comment add all sets for the current interval and continue search
\State add all $\setu$ to $\setU^\mathrm{opt}(\e)$
\State $T \leftarrow T - T_j$
\EndIf
\EndFor
\end{algorithmic}
\end{algorithm}

\begin{algorithm}[!h]
\caption{(Subroutine: \textbf{opt-search})}\label{alg:search-opt}
\textbf{inputs}: $\setU^\mathrm{opt}(\e)$, $T_j$, $(\bar{\gamma}_\setu)_{|\setu| < \infty}$, 
$L_j^\mathrm{unsorted}$, $L_{j + 1}$, $I_j$\\
\textbf{outputs}: $\ell_\mathrm{next}$, $T_j$, $L_j^\mathrm{unsorted}$, $L_{j+ 1}$
\begin{algorithmic}[1]
\State $\ell_\mathrm{next} = 1$
\For{{$\setu \in L_j$}}
\State $i \leftarrow |\setu|$
\While{{$i > 0$}}
\State $L_j \leftarrow L_j \setminus \setu$
\Comment reducing the double-handling of sets
\If{{$\bar{\gamma}_\setu \in I_j$}}
\State \textbf{if} $\setu \in L_j^\mathrm{unsorted}$ \textbf{then} break
\Comment already visited $\setu$ and any future increments
\State add $\setu$ to $L_j^\mathrm{unsorted}$
\State $T_j \leftarrow T_j + \bar{\gamma}_\setu$
\State $i = |\setu|$
\Else
\State add $\setu$ to $L_{j + 1}$
\Comment $\setu$ to be checked first in next interval
\State $i \leftarrow i - 1$
\State \textbf{if} $i = 0$ \textbf{then} break
\Comment go to next $\setu \in L_j$
\EndIf
\State $\setu \leftarrow $\textbf{increment-$\setu$}$(\setu, i)$
\EndWhile
\State $\ell_\mathrm{next} = |\setu| + 1$
\Comment main search will start at cardinality $|\setu| + 1$
\EndFor
\State return $(\ell_\mathrm{next}, T_j, L_j^\mathrm{unsorted}, L_{j + 1})$
\end{algorithmic}
\end{algorithm}

\section{Discussion}
In this paper we have introduced the notion of \emph{superposition dimension}
and \emph{optimal} active sets to be used in the MDM
for multivariate integration and presented an algorithm detailing
their construction.
We also introduced a second simplified, computationally less intensive 
version of the algorithm, which constructs \textit{quasi-optimal} active sets.
Our numerical results show that the quasi-optimal active sets are of a
similar size to the optimal active sets. Often the two sets are exactly the same.
In all of our numerical results the optimal and quasi-optimal active sets
are smaller than, and have superposition dimension less than or equal to,
the active sets using the construction in \cite{PlaWas10}.

To observe how
different choices of parameters $a$ and $c$
affect our construction, statistics 
 on the resulting optimal active sets are given in Tables~\ref{tab:U_p=2}-\ref{tab:d_p=inf}. Tables~\ref{tab:U_p=2} and
\ref{tab:d_p=2} give, respectively, the size and the superposition dimension
of the optimal active set for $p = 2$ and an error request of $10^{-2}$.
For $p = \infty$, $\e = 10^{-2}$ the size and superposition dimension
of the optimal active sets are given in Tables~\ref{tab:U_p=inf} and \ref{tab:d_p=inf}.
The results for the quasi-optimal active set are again very similar and so
have not been included here.
As expected these results demonstrate that as
the decay of the weights is slower or the weights become larger  
($a$ smaller and $c$ larger)
the problem becomes more difficult and the active sets are by necessity larger.
However the superposition dimension remains relatively small, 
at most~6.

\begin{table}[!h]
\centering
\begin{minipage}{.4\textwidth}
\centering
\begin{tabular}{l||c|c|c}
 & \multicolumn{3}{c}{$a$}\\
\hline
$c$ & 4 & 3 & 2\\
\hline
$\tfrac{1}{2}$ & 3 & 5 & 12
\\
1 & 4 & 7 & 30
\\
2 & 6 & 14 & 122
\end{tabular}
\caption{$|\setU^\mathrm{opt}(10^{-2})|$ for $p = 2$
\hfill
and different $a$, $c$.}\label{tab:U_p=2}
\end{minipage}
\hspace{0.05\textwidth}
\begin{minipage}{.4\textwidth}
\centering
\begin{tabular}{l||c|c|c}
 & \multicolumn{3}{c}{$a$}\\
\hline
$c$ & 4 & 3 & 2\\
\hline
$\tfrac{1}{2}$ & 1 & 2 & 2
\\
1 & 2 & 2 & 3
\\
2 & 2 & 3 & 4
\\
\end{tabular}
\caption{$d(\setU^\mathrm{opt}(10^{-2}))$ for $p = 2$
and different $a$, $c$.}\label{tab:d_p=2}
\end{minipage}
\end{table}

\begin{table}[!h]
\centering
\begin{minipage}{.4\textwidth}
\centering
\begin{tabular}{l||c|c|c}
 & \multicolumn{3}{c}{$a$}\\
\hline
$c$ & 4 & 3 & 2\\
\hline
$\tfrac{1}{2}$ & 4 & 7 & 150
\\
1 & 5 & 15 & 1346
\\
2 & 8 & 43 & 31,013
\end{tabular}
\caption{$|\setU^\mathrm{opt}(10^{-2})|$ for $p = \infty$
and different $a$, $c$.}\label{tab:U_p=inf}
\end{minipage}
\hspace{.05\textwidth}
\begin{minipage}{.4\textwidth}
\centering
\begin{tabular}{l||c|c|c}
 & \multicolumn{3}{c}{$a$}\\
\hline
$c$ & 4 & 3 & 2\\
\hline
$\tfrac{1}{2}$ & 2 & 3 & 4
\\
1 & 2 & 2 & 4
\\
2 & 2 & 3 & 6
\end{tabular}
\caption{$d(\setU^\mathrm{opt}(10^{-2}))$ for $p = \infty$
and different $a$, $c$.}\label{tab:d_p=inf}
\end{minipage}
\end{table}
\pagebreak
Finally, we have constructed the sets $\setU^{\rm q-opt}(\e)$,
$\setU^{\rm opt}(\e)$, and $\setU^{\rm PW}(\e)$  for the
product weights with $c=1$. They are listed explicitly in the Appendix. 

\section*{Acknowledgements}
The authors would like to thank the two anonymous referees
for their constructive comments which helped improve the paper.

\section*{Appendix}
We list here a selection of the constructed active sets 
from the previous sections (the very largest sets have been omitted).
To save the space sometimes we write
$[...\{x_1,\dots,x_k,x_{k+1}\}]$ to denote the sequence of sets,
\[
  [...\{x_1,\dots,x_k,x_{k+1}\}]\,=\,\{x_1,\dots,x_k,x_k+1\},\dots,
  \{x_1,\dots,x_k,x_{k+1}\}. 
\]
For instance $[...\{3\}]$ denotes $\{1\},\{2\},\{3\}$ and 
$[...\{1,5\}]$ denotes $\{1,2\},\{1,3\},\{1,4\},\{1,5\}$.

\medskip\noindent
{\bf Case $p=1$ and $a=4$}
{\footnotesize
\[
\setU^{\rm PW}(10^{-1})\,=\,\{\emptyset,\{1\}\},
\quad
\setU^{\rm PW}(10^{-2})\,=\,\{\emptyset,[...\{3\}],\{1,2\},\{1,3\}\},
\quad 
\setU^{\rm PW}(10^{-3})\,=\,\{\emptyset,[...\{5\}],[...\{1,5\}]\}.
\]
}
{\bf Case of $p=1$ and $a=3$}
{\footnotesize
\[
\setU^{\rm PW}(10^{-1})\,=\,\{\emptyset,\{1\},\{2\},\{1,2\}\},
\quad 
\setU^{\rm PW}(10^{-2})\,=\,\{\emptyset,[...\{4\}],
[..\{1,4\}]\}, 
\]
\[
\setU^{\rm PW}(10^{-3})\,=\,\{\emptyset,[...\{9\}],[...\{1,9\}],
\{2,3\},\{2,4\},\{1,2,3\},\{1,2,4\}\}.
\]
}
{\bf Case $p=1$ and $a=2$}
{\footnotesize
\[
\setU^{\rm PW}(10^{-1})\,=\,\{\emptyset,[...\{3\}],
\{1,2\},\{1,3\}\},
\]
\[
\setU^{\rm PW}(10^{-2})\,=\,\{\emptyset,[...\{9\}],
[...\{1,9\}],\{2,3\},\{2,4\},\{1,2,3\},\{1,2,4\},
\]
\begin{eqnarray*}
\setU^{\rm PW}(10^{-3})&=&\{\emptyset,[...\{31\}],
[...\{1,31\}],[...\{2,15\}],[...\{3,10\}],[...\{4,7\}],\{5,6\},
[...\{1,2,15\}],\\
  &&\quad[...\{1,3,10\}],[...\{1,4,7\}]
  ,\{1,5,6\},\{2,3,4\},\{2,3,5\}, \{1, 2, 3, 4\}, \{1, 2, 3, 5\}\}.
\end{eqnarray*}
}
{\bf Case $p=2$ and $a=4$}
{\footnotesize 
\[
\setU^{\rm q-opt}(10^{-1})\,=\,\setU^{\rm opt}(10^{-1})
\,=\,\{\emptyset,\{1\}\},
\quad\mbox{and}\quad
\setU^{\rm PW}(10^{-1})\,=\,\{\emptyset,\{1\},\{2\}\},
\]
\[
\setU^{\rm q-opt}(10^{-2})\,=\,\setU^{\rm opt}(10^{-2})\,=\,
\{\emptyset,\{1\},\{2\},\{1,2\}\}\quad\mbox{and}\quad
\setU^{\rm PW}(10^{-2})\,=\,\{\emptyset,[...\{4\}],[...\{1,4\}]\},
\]
\[
\setU^{\rm q-opt}(10^{-3})\,=\,\setU^{\rm opt}(10^{-3})\,=\,
\{\emptyset,[...\{5\}],[...\{1,4\}]\}
\]
\[
\setU^{\rm PW}(10^{-3})\,=\,\{\emptyset,[...\{9\}],[...\{1,8\}],
\{2,3\},\{2,4\}, \{1, 2, 3\}\}.
\]
}
{\bf Case of $p=2$ and $a=3$}
{\footnotesize
  \[
\setU^{\rm q-opt}(10^{-1})\,=\,\setU^{\rm opt}(10^{-1})\,=\,
\{\emptyset,\{1\}\}\quad\mbox{and}\quad
\setU^{\rm PW}(10^{-1})\,=\,\{\emptyset,[...\{3\}],\{1,2\}\},
\]
\[
  \setU^{\rm q-opt}(10^{-2})\,=\,\setU^{\rm opt}(10^{-2})\,=\,
\{\emptyset,[...\{4\}],\{1,2\},\{1,3\}\},
\]
\[
\setU^{\rm PW}(10^{-2})\,=\,
\{\emptyset,[...\{9\}],[...\{1,7\}],\{2,3\},\{1,2,3\}\},
\]
\[
\setU^{\rm q-opt}(10^{-3})\,=\,\{\emptyset,[...\{12\}],[...\{1,10\}],
[...\{2,5\}],\{1,2,3\}\},
\]
\[
  \setU^{\rm opt}(10^{-3})\,=\,\{\emptyset,[...\{11\}],[...\{1,9\}],
  \{2,3\},\{2,4\},\{1,2,3\},\{1,2,4\}\},
\]
\[
\setU^{\rm PW}(10^{-3})\,=\,\{\emptyset,[...\{26\}],[...\{1,21\}],
[...\{2,10\}],[...\{3,7\}],\{4,5\},[...\{1,2,9\}],[...\{1,3,6\}]\}.
\]
}
{\bf Case of $p=2$ and $a=2$}
{\footnotesize
\[
\setU^{\rm opt}(10^{-1})\,=\,\{\emptyset,\{1\},\{2\},\{1,2\}\}
\quad\mbox{and}\quad   
\setU^{\rm q-opt}(10^{-1})\,=\,\{\emptyset,[...\{4\}],\{1,2\}\},
\]
\[\setU^{\rm PW}(10^{-1})\,=\,\{\emptyset,[...\{8\}],[...\{1,6\}], \{2, 3\}\}
\]
\[
\setU^{\rm q-opt}(10^{-2})\,=\,
\{\emptyset,[...\{18\}],[...\{1,10\}],[...\{2,5\}], \{1, 2, 3\}\},
\]
\[
\setU^{\rm opt}(10^{-2})\,=\,\{\emptyset,[...\{14\}],[...\{1,11\}],
[...\{2,5\}],[...\{1,2,4\}]\},
\]
\begin{eqnarray*}
\setU^{\rm PW}(10^{-2})&=&\{\emptyset,[...\{54\}],[...\{1,41\}],
[...\{2,20\}],[...\{3,13\}],[...\{4,10\}],[...\{5,8\}],
,[...\{1,2,15\}],\\
&&\quad[...\{1,3,10\}],[...\{1,4,7\}],\{1,5,6\},\{2,3,4\},\{2,3,5\}\}.
\end{eqnarray*}
}
{\bf Case $p=\infty$ and $a=4$}
{\footnotesize
\[
\setU^{\rm q-opt}(10^{-1})\,=\,\setU^{\rm opt}(10^{-1})\,=\,
\{\emptyset,\{1\}\} \quad\mbox{and}\quad
\setU^{\rm PW}(10^{-1})\,=\,\{\emptyset,[...\{4\}],
\{1,2\},\{1,3\}\}.
\]
\[
\setU^{\rm q-opt}(10^{-2})\,=\,\setU^{\rm opt}(10^{-2})\,=\,
\{\emptyset,[...\{3\}],\{1,2\}\},
\]
\[\setU^{\rm PW}(10^{-2})\,=\,\{\emptyset,[...\{10\}],
[...\{1,8\}],\{2,3\},\{2,4\},\{1,2,3\}\},
\]
\[
\setU^{\rm q-opt}(10^{-3})\,=\,
\{\emptyset,[...\{8\}],[...\{1,7\}]\},
\]
\[
\setU^{\rm opt}(10^{-3})\,=\,
\{ \emptyset, [...\{8\}], [...\{1, 6\}], \{2, 3\}\}
\]
\[
  \setU^{\rm PW}(10^{-3})\,=\,\{\emptyset,[...\{26\}],
[...\{1,22\}],[...\{2,11\}],[...\{3,7\}],\{4,5\},[...\{1,2,9\}],
[...\{1,3,6\}]\}
\]
}
{\bf Case $p=\infty$ and $a=3$}
{\footnotesize
\[\setU^{\rm q-opt}(10^{-1})\,=\,\setU^{\rm opt}(10^{-1})\,=\,
\{\emptyset,\{1\},\{2\}\},
\]
\[\setU^{\rm PW}(10^{-1})\,=\,\{\emptyset,[...\{10\}],[...\{1,8\}],
\{2,3\},\{2,4\},\{1,2,3\}\},
\]
\[
\setU^{\rm q-opt}(10^{-2})\,=\,\setU^{\rm opt}(10^{-2})\,=\,
\{\emptyset,[...\{8\}],[...\{1,6\}],\{2,3\}\},
\]
\begin{eqnarray*}
&&\setU^{\rm PW}(10^{-2}) \,=\, \{\emptyset, 
[...\{49\}], [...\{1, 39\}], [...\{2, 19\}], [... \{3, 13\}], [...\{4, 9\}],
\{5, 6\},\\
&&\{5, 7\}, [...\{1, 2, 15\}], [...\{1, 3, 10\}], [...\{1, 4, 7\}],
\{1, 5, 6\}, \{2, 3, 4\}, \{2, 3, 5\}, \{1, 2, 3, 4\}
\},
\end{eqnarray*}
\[
\setU^{\rm q-opt}(10^{-3})\,=\,
\{\emptyset,[...\{36\}],[...\{1,29\}],[...\{2,14\}]
,[...\{3, 9\}], [...\{4, 7\}], 
[...\{1,2,8\}]\},
\]
\[
\setU^{\rm opt}(10^{-3})\,=\,\{\emptyset,[...\{31\}],[...\{1,25\}],
[...\{2,12\}],[...\{3, 8\}], [...\{4, 6\}],
[...\{1,2,9\}],[...\{1,3,6\}]\},
\]}
\begin{eqnarray*}
&&\setU^{\rm PW}(10^{-3})\,=\,\{\emptyset,[...\{208\}],
[...\{1,165\}],[...\{2,82\}],[...\{3,55\}],[...\{4,41\}],
[...\{5,33\}],[...\{6,27\}],\\
&&[...\{7,23\}],[...\{8,20\}],[...\{9,18\}],[...\{10,16\}],
[...\{11,15\}],\{12,13\},[...\{1,2,65\}],[...\{1,3,43\}],\\
&&[...\{1,4,32\}],[...\{1,5,26\}],
[...\{1,6,21\}],[...\{1,7,18\}],[...\{1,8,16\}],
[...\{1,9,14\}],[...\{1,10,13\}],\\
&&[...\{2,3,21\}],
[...\{2,4,16\}],[...\{2,5,13\}],[...\{2,6,10\}],
\{2,7,8\},\{2,7,9\},[...\{3,4,10\}],[...\{3,5,8\}],\\
&&\{3,6,7\},\{4,5,6\},[...\{1,2,3,17\}],
[...\{1,2,4,13\}],[...\{1,2,5,10\}],\{1,2,6,7\},\{1,2,6,8\},\\
&&[...\{1,3,4,8\}],\{1,3,5,6\}\}
\end{eqnarray*}
{\bf Case $p=\infty$ and $a=2$}
{\footnotesize 
\[
\setU^{\rm q-opt}(10^{-1})\,=\,\{\emptyset,[...\{22\}],
[...\{1,15\}],\{2,3\}\}
\]
\[
\setU^{\rm opt}(10^{-1})\,=\,\{\emptyset,[...\{16\}],
[...\{1,12\}],[...\{2,5\}], [...\{1, 2, 4\}]\},
\]
\begin{eqnarray*}
&&\setU^{\rm PW}(10^{-1}) \,=\,
\{\emptyset, [...\{511\}], [...\{1, 361\}], [..., \{2, 180\}], [...\{3, 120\}], [...\{4, 90\}], [...\{5, 72\}], [...\{6, 60\}],\\
&&[...\{7, 51\}], [...\{8, 45\}], [...\{9, 40\}], [...\{10, 36\}], [...\{11, 36\}], [...\{11, 32\}], [...\{12, 30 \}],[...\{ 13, 27\}],\\
&& [...\{ 14, 25\}], [...\{ 15, 24\}], [...\{ 16, 22\}],[...\{17, 21\}], [...\{ 18, 20\}],[...\{ 1, 2, 127\}], [...\{1, 3, 85\}],\\
&& [...\{1, 4, 63\}],[...\{1, 5, 51\}], [...\{1, 6, 42\}], [...\{1, 7, 36\}],[...\{1, 8, 31\}], [...\{1, 9, 28\}],[...\{1, 10, 25\}],\\
&&[...\{1, 11, 23\}], [...\{1, 12, 21\}], [...\{1, 13, 19\}],[...\{1, 14, 18\}],[...\{1, 15, 17\}], [...\{2, 3, 42\}],\\
&&[...\{2, 4, 31\}], [...\{2, 5, 25\}], [...\{2, 6, 21\}],[...\{2, 7, 18\}], [...\{2, 8, 15\}], [...\{2, 9, 14\}], [...\{2, 10, 12\}],\\
&&[...\{3, 4, 21\}],[...\{3, 5, 17\}], [...\{3, 6, 14\}], [...\{3, 7, 12\}], [...\{3, 8, 10\}], [...\{4, 5, 12\}],[...\{4, 6, 10\}],\\
&&[...\{4, 7, 9\}], [...\{5, 6, 8\}], [...\{1, 2, 3, 30\}], [...\{1, 2, 4, 22\}],[...\{1, 2, 5, 18\}], [...\{1, 2, 6, 15\}],\\
&&[...\{1, 2, 7, 12\}], [...\{1, 2, 8, 11\}], \{1, 2, 9, 10\}, [...\{1, 3, 4, 15\}], [...\{1, 3, 5, 12\}], [...\{1, 3, 6, 10\}],\\
&&\{1, 3, 7, 8\}, [...\{1, 4, 5, 9\}],\{1, 4, 6, 7\}, [...\{2, 3, 4, 7\}], \{2, 3, 5, 6\}, \{1, 2, 3, 4, 5\}\}.
\end{eqnarray*}
}

%------------------------------------------------------------------------------
%------------------------------------------------------------------------------

\begin{small}
\noindent{\bf Authors' Addresses:}\\
A.~D.~Gilbert\\
School of Mathematics and Statistics\\
The University of New South Wales\\
Sydney, NSW 2052, Australia\\
E-mail: \texttt{alexander.gilbert@student.unsw.edu.au}\\

\noindent G.~W.~Wasilkowski\\
Department of Computer Science\\
 University of Kentucky\\
301 David Marksbury Building\\
 Lexington, KY 40506, USA\\
E-mail: \texttt{greg@cs.uky.edu}\\

\noindent {\bf Funding:}
The author A. D. Gilbert is grateful for the financial support received from the Australian Research Council (DP150101770) and the University of New South Wales.
\end{small}
\end{document}